\documentclass[12pt,emtex]{report}
\usepackage{times}

\usepackage{graphicx}
\usepackage{amsfonts}       %
\usepackage{amsmath} 
\usepackage{amssymb}        %
\usepackage{amsthm}         

\usepackage{epsfig}
\psfigdriver{dvips}

\usepackage{latexsym}
\usepackage[matrix,curve]{xypic}

\usepackage[german]{babel}  
\newcommand{\Proj}{\mathbb{P}}

\newcommand{\Z}{\mathbb{Z}}
\newcommand{\Q}{\mathbb{Q}}

\DeclareMathSizes{12}{11}{8}{5}

\begin{document}

\title{Torelli's theorem from the topological point of view}
\maketitle

\bigskip
Torelli's theorem states, that the isomorphism class of a smooth
projective curve of genus $g\geq 2$ over an algebraically closed
field $k$ is uniquely determined by the isomorphism class of the
associated pair $(X,\Theta)$, where $X$ is the Jacobian variety of
$C$ and $\Theta$ is the canonical theta divisor. The aim of this
note is to give a \lq{topological\rq}\ proof of this theorem.
Although Torelli's theorem is not a topological statement the
proof to be presented gives a characterization of $C$ in terms of
perverse sheaves on the Jacobian variety $X$, which are attached
to the theta divisor by a \lq{topological}\rq\ construction.

\bigskip
For complexes $K,L\in D_c^b(X,\overline\Q_l)$ define $K*L\in D_c^b(X,\overline\Q_l)$ by the
direct image complex $Ra_*(K\boxtimes L)$, where $a:X\times X\to X$ is the addition law of $X$.
Let $K^0_*(X)$ be the tensor product of the Grothendieck group of perverse sheaves on $X$ with
the polynomial ring $\Z[t^{1/2},t^{-1/2}]$. $K^0_*(X)$ is a commutative ring with ring
structure defined by the convolution product, hence also the  quotient ring $K_*(X)$  obtained
by dividing the principal ideal generated by the constant perverse sheaf $\delta_X$ on $X$.
Both rings $K_*^0(X)$ and $K_*(X)$ resemble properties of the homology ring of $X$ endowed with
the $*$-product, but have a much richer structure. A sheaf complex $L\in
D_c^b(X,\overline\Q_l)$ defines a class in $K_*(X)$ by the perverse Euler characteristic
$\sum_\nu (-1)^\nu\cdot {}^p\! H^\nu(L)\cdot t^{\nu/2}$ . Similar to the homology ring every
irreducible closed subvariety $Y$ has a class in $K_*(X)$ defined as the class of the perverse
intersection cohomology sheaf $\delta_Y$ of $Y$. For details see [W]. This allows to consider
the product
$$ \delta_\Theta * \delta_\Theta \in K_*(X) \ .$$

Whereas the corresponding product in the homology ring of $X$ is
zero, this product turns out to be nonzero. $\delta_\Theta
* \delta_\Theta$ is of the form $\sum_{\nu,\mu} A_{\nu,\mu}
t^{\mu/2}$. Recall, the coefficients are irreducible perverse
sheaves $A_{\nu,\mu}$ on $X$. For a perverse sheaf $A$ on $X$,
which  is a sheaf complex on $X$,  let ${\cal H}^i(A)$ denote the
associated cohomology sheaves for $i\in \Z$. Let $\kappa\in X(k)$
be the Riemann constant defined by $\Theta = \kappa - \Theta$. It
depends on the choice of the Abel-Jacobi map $C\to X$.

\bigskip
{\bf Theorem}: {\it Let $C$ be a curve of genus $g\geq 3$. There
exists a unique irreducible perverse sheaf $A=A_{\nu,0}$, among
the coefficients of $\delta_\Theta*\delta_\Theta$, characterized
by one of the following equivalent properties
\begin{enumerate} \item ${\cal H}^{-1}(A)$ is nonzero, but
not a constant sheaf on $X$.
\item ${\cal H}^{-1}(A)$  is the skyscraper sheaf $H^1(C)\otimes \delta_{\{\kappa\}}$ with
support in the point $\kappa\in X$.
\end{enumerate}
Furthermore the support of the perverse sheaf $A$ is
$\kappa+C-C\subseteq X$.}

\bigskip
Taking this for granted, Torelli's theorem is an immediate
consequence. In fact the subvariety $C-C$ uniquely determines the
curve $C$. This is well known. For instance, $C$ is the unique
one-dimensional fiber of the minimal desingularization of $C-C$.
Since the theorem stated above allows to recover the sheaf complex
$A$ from $\Theta$, this determines $C$  from the data $(X,\Theta)$
via the support $\kappa+C-C$ of $A$ for $g\geq 3$. The cases
$g=1,2$ are trivial.

\bigskip
\underbar{Remark}: $A$ is a direct summand of the complex
$\delta_{C}*\delta_{\kappa -C}$. If $C$ is not hyperelliptic, it
splits into two irreducible perverse summands $\delta_{\{\kappa\}}
\oplus A$. If $C$ is hyperelliptic, then it splits into the three
irreducible perverse summands $\delta_{\{\kappa\}} \oplus
\delta_{\kappa+C-C}\oplus A$.

\bigskip
We now give a sketch of the theorem in the non-hyperelliptic case.
For all integers $r\geq 0$ let $\delta_r\in K_*(X)$ be the class
of the direct image complex $Rp_{r,*}\delta_{C^{(r)}}$, where
$p_r:C^{(r)}=C^r/\Sigma_r \to X$ are the higher Abel-Jacobi maps
from the symmetric quotient of $C^r$ to $X$. For $r\leq g-1$ the
image of $p_r$ is the Brill-Noether subvariety $W_r=C+\cdots + C$
($r$ copies) of $X$. If $C$ is not hyperelliptic, then
$$ \delta_{W_r} =\delta_r \ ,$$ since $p_r$ is a small morphism
by the theorem of Martens [M] for $r\leq g-1$. In particular
$\delta_\Theta=\delta_{g-1}$, which will be used in the proof. (In
the hyperelliptic case $\delta_\Theta=\delta_{g-1}-\delta_{g-3}$.
For this and further details we refer to [W]).

\bigskip
\underbar{Proof of the theorem}: Suppose $C$ is not hyperelliptic.

\bigskip
1) Since the canonical morphism
$$\tau:C^{(i)}\times C^{(j)}\to C^{(i+j)}$$ is a finite ramified
covering map, the direct image $R\tau_* \delta_{C^{(i)}\times C^{(j)}}$ decomposes into a
direct sum of etale sheaves $ \bigoplus_{\nu} m(i,j,\nu)\cdot {\cal F}_{i+j-\nu,\nu}$ by
keeping track of the underlying action of the symmetric group $\Sigma_{i+j}$ for the map
$C^{i+j}\to C^{(i+j)}$ (see [W].4.1). If we apply $Rp_{i+j,*}$, this gives a formula for
$\delta_i* \delta_j$. From $p_{i+j}\circ \tau = a\circ (p_i\times p_j)$, where $a:X\times X\to
X$ is the addition law of $X$, one obtains for $i\geq j$ that the convolution $\delta_i
* \delta_j$ is $ \delta_{i+j} \oplus \delta_{i+j-1,1} \oplus \cdots
\oplus \delta_{i-j,j} $, where $\delta_{r,s}= Rp_{i+j,*}({\cal
F}_{r,s})$. A special case is
$$    \delta_\Theta * \delta_\Theta= \delta_{g-1}*\delta_{g-1} = \delta_{2g-2} \oplus \delta_{2g-3,1}
\oplus \cdots \oplus \delta_{g-1,g-1}  \ .$$ Another case is $ \
\delta_1
* \delta_{2g-3} = \delta_{2g-2} \oplus \delta_{2g-3,1}$, and
together this implies $$ \fbox{$
\delta_{2g-3}*\delta_1\hookrightarrow \delta_\Theta*\delta_\Theta
$} \ .$$

\bigskip
2) The morphism $f:C\times C\to \kappa+C-C \subseteq X$, defined
by $(x,y)\mapsto \kappa+x-y$, is semi-small. If $C$ is not
hyperelliptic, then $f$ is a birational map, which blows down the
diagonal to the point $\kappa$, and is an isomorphism otherwise.
Hence the direct image $Rf_*(\delta_C\boxtimes\delta_C)$  is
perverse on $X$, and necessarily decomposes
$Rf_*(\delta_C\boxtimes\delta_C) =\delta_C* \delta_{\kappa-C} =
\delta_{\{\kappa\}} \oplus \delta_{\kappa+C-C}$ such that
$$ \fbox{$ {\cal H}^{-1}(\delta_{\kappa+C-C}) \cong H^1(C)\otimes \delta_{\{\kappa\}} \
$} \ .$$

\bigskip

3) We claim $\delta_{2g-3} \equiv \delta_{\kappa-C}$ and $\delta_{2g-2} \equiv
\delta_{\{\kappa\}}$ in $K_*(X) $ (ignoring Tate twists). These are the simplest cases of the
duality theorem [W] 5.3. This implies
$$ \fbox{$ \delta_{2g-3}*\delta_1 \ \equiv\ \ \ \delta_{\{\kappa\}} + \delta_{\kappa+C-C} $} \ ,$$
in $K_*(X)$ using step 2.

\bigskip
\underbar{Proof of the claim}: By the theorem of Riemann-Roch $
C^{(2g-3)} \overset{p}{\to} X $ is a $\Proj^{g-2}$-bundle over
$\kappa-C$ and a $\Proj^{g-3}$-bundle over the open complement
$X\setminus \ (\kappa-C)$. Hence $Rp_*\delta_{C^{(2g-3)}}$ is a
direct sum of $ \delta_{\kappa - C}$ and a sum of translates of
constant sheaves on $X$.  Similarly
$pr_{2g-2}^{-1}(\{\kappa\})=\Proj^{g-1}$, and $pr_{2g-2}$ is a
$\Proj^{g-2}$-bundle over the open complement $X\setminus \
\{\kappa\}$. Hence $\delta_{2g-2}\equiv \delta_{\{\kappa\}} $ in
$K_*(X)$.

\bigskip
4) $\Theta=\kappa-\Theta$ and the definition of the convolution
product
 implies
$$ {\cal H}^{-1}(\delta_\Theta * \delta_\Theta) \ \cong \
IH^{2g-1}(\Theta)\otimes \delta_{\{\kappa\}} \ $$ for an arbitrary principally polarized
abelian varieties $(X,\Theta)$, where $IH^{2g-1}(\Theta)$ denotes the intersection cohomology
group of $\Theta$. If the singularities of $\Theta$ have codimension $\geq 3$ in $\Theta$, then
ignoring Tate twists this implies (see [W] {2.9} and [W2])
$$ \fbox{$ {\cal H}^{-1}(\delta_\Theta * \delta_\Theta) \ \cong \
H^1(X)\otimes \delta_{\{\kappa\}} $} \ .$$

\bigskip In fact
by the Hard Lefschetz theorem $IH^{2g-3}(\Theta)$ and $H^1(X)$
have the same dimensions. A more elementary argument proves this
for Jacobians including the case of hyperelliptic curves: For
example  for non-hyperelliptic curves we have $IH^{\bullet}(W_d) =
H^{\bullet}(X,Rp_{d,*}\delta_{C^{(d)}}[-d])=H^{\bullet}(C^{(d)})=(\bigotimes^d
H^\bullet(C))^{\Sigma_d}$, since $p_{d}$ is a small morphism. Thus
$$IH^{d+\bullet}(W_d)\ \cong\ \bigoplus_{a+b=d}\ Sym^a\Bigl(H^0(C)[1]\oplus
H^2(C)[-1]\Bigr) \otimes \Lambda^b(H^1(C))\ .$$ For $IH^{2d-1}(W_d)$ only $a=d-1$ contributes,
hence $IH^{2d-1}(W_d)\!\cong\! H^1(C)\!\cong\! H^1(X)$. (For the hyperelliptic case see [W]
4.2).

\bigskip
\underbar{Conclusion}: For curves $C$, which are not
hyperelliptic, the perverse sheaf $A$ defined by
$\delta_{\kappa+C-C}$ satisfies all the assertions of the theorem.
$\delta_{\kappa+C-C}$ is a direct summand of
$\delta_\Theta*\delta_\Theta$ by step 1 and 3. By step 2 and 4 we
obtain modulo constant sheaves on $X$
$${\cal H}^{-1}(\delta_\Theta*\delta_\Theta) \equiv H^1(X)\otimes
\delta_{\{\kappa\}} \equiv H^1(C)\otimes \delta_{\{\kappa\}}
\equiv {\cal H}^{-1}(\delta_{\kappa+C-C}) \ .$$ Since $K_*(X)$ is
a quotient of $K_*^0(X)$, the last identity only holds modulo
constant sheaves on $X$. But this suffices to imply the theorem.

\bigskip
\underbar{Remark}: In [W] we constructed a $\overline\Q_l$-linear
Tannakian category ${\cal BN}$ attached to $C$  equivalent to the
category of finite dimensional $\overline\Q_l$-representations
$Rep(G)$, where $G$ is $Sp(2g-2,\overline\Q_l)$ or
$Sl(2g-2,\overline\Q_l)$ depending on whether $C$ is hyperelliptic
or not. In this category $\delta_\Theta$ corresponds to the
alternating power $\Lambda^{g-1}(st)$ of the standard
representation, and $A$ corresponds to the adjoint representation.

\bigskip\noindent
{\bf Bibliography}

\bigskip\noindent
[M] Martens H.H, On the variety of special divisors on a curve,
Crelle 227 (1967), p.111 -- 120.

\bigskip\noindent
[W] Weissauer R., Brill-Noether sheaves (preprint)

\bigskip\noindent
[W2] Weissauer R., Inner cohomology (preprint)

\end{document}